\RequirePackage{fix-cm}
\documentclass[smallcondensed]{svjour3}     
\usepackage{graphicx}%
\usepackage{multirow}%
\usepackage{amsmath,amssymb,amsfonts}%
\usepackage{mathrsfs}%
\usepackage[title]{appendix}%
\usepackage{xcolor}%
\usepackage{textcomp}%
\usepackage{manyfoot}%
\usepackage{booktabs}%
\usepackage{algorithm}%
\usepackage{algorithmicx}%
\usepackage{algpseudocode}%
\usepackage{listings}%

\usepackage{mathtools}

\usepackage{tikz}
\usetikzlibrary{calc}
\usepackage{pgfplots}
\usepgfplotslibrary{fillbetween}
\pgfplotsset{compat=1.18} 

\usepackage{subfig}
\usepackage{float}
\usepackage{adjustbox}
\usepackage{hyperref}
\usepackage[numbers]{natbib}

\usepackage{nomencl}
\usepackage{framed}
\usepackage{etoolbox} 
\usepackage{multicol}
\makenomenclature
\patchcmd{\thenomenclature}{\section*{\nomname}}{}{}{}
\renewcommand{\nomname}{}          

\renewcommand\nomgroup[1]{%
  \item[\bfseries
  \ifstrequal{#1}{A}{Sets}{%
  \ifstrequal{#1}{B}{Constants}{%
  \ifstrequal{#1}{C}{Variables}{%
  \ifstrequal{#1}{D}{Functions}{}}}}%
  ]}

\def\eg{\emph{e.g., }}
 
\def\ea{\emph{et al. }} 

\DeclareMathOperator*{\subjto}{subject\,to}

\begin{document}

\title{A Quadratically-Constrained Convex Approximation for the AC Optimal Power Flow\thanks{
The authors would like to acknowledge the support from Davidson School of Chemical Engineering at Purdue University and the Office of Naval Research (grant no. N000142412641).
}
}

\author{Gonzalo Constante-Flores \and
        Can Li}


\institute{Gonzalo Constante-Flores \and Can Li \at
              Davidson School of Chemical Engineering, Purdue University,\\
              610 Purdue Mall, West Lafayette, IN 47907, USA\\
              \email{\{geconsta,canli\}@purdue.edu}
}

\date{Received: date / Accepted: date}

\maketitle

\begin{abstract}
We introduce a quadratically-constrained approximation (QCAC) of the AC optimal power flow (AC-OPF) problem. Unlike existing approximations like the DC-OPF, our model does not rely on typical assumptions such as high reactance-to-resistance ratio, near-nominal voltage magnitudes, or small angle differences, and preserves the structural sparsity of the original AC power flow equations, making it suitable for decentralized power systems optimization problems. To achieve this, we reformulate the AC-OPF problem as a quadratically constrained quadratic program. The nonconvex terms are expressed as differences of convex functions, which are then convexified around a base point derived from a warm start of the nodal voltages. If this linearization results in a non-empty constraint set, the convexified constraints form an inner convex approximation. Our experimental results, based on Power Grid Library instances of up to 30,000 buses, demonstrate the effectiveness of the QCAC approximation with respect to other well-documented conic relaxations and a Taylor-series linear approximation. We further showcase its potential advantages over the well-documented second-order conic relaxation of the power flow equations in the optimal reactive power dispatch and photovoltaic hosting capacity problem.
\subclass{90C26 \and 90C90 \and 65K05}
\end{abstract}

\section{Introduction} \label{sec:introduction}

Power systems are undergoing important changes due to the pervasive electrification of transportation and industry sectors and the increasing share of renewable energy resources. System operators must rely on reliably solving the AC Optimal Power Flow (AC-OPF) problem to ensure economically efficient and secure planning and operation decisions. The main challenge of reliably solving the AC-OPF comes from the nonconvexity of the power flow equations, which describe the steady-state governing physical laws of the transmission network \cite{Capitanescu2016}.

In recent years, significant efforts have been devoted to developing more tractable linear and convex relaxations of the power flow equations, which can be reliably solved using state-of-the-art optimization solvers. In particular, these efforts include the semi-definite programming \cite{Bai2008}, second-order cone \cite{Jabr2006}, and the quadratic convex \cite{Coffrin2016} relaxations. The key features of convex relaxations are their ability to provide infeasibility certificates, lower bounds on the optimal objective values, and globally optimal decisions when tight. The main drawbacks of using relaxations include the overestimation of power losses \cite{Chen2016}, relaxations of Ohm's law \cite{Sojoudi2012}, and their inability to leverage AC base solutions, also known as warm starts, from previously solved instances. 

Another line of research is focused on developing tractable model-based approximations of the power flow equations \cite{Molzahn2019}. One of the main advantages of such approximations is that the can (i) use solutions of previous operating conditions \cite{Coffrin2014}, (ii) leverage variable estimates from the solution of the state estimator for real-time applications, (iii) employ engineering assumptions \cite{Alsac1974}, \eg reactance-to-resistance ratio, flat voltage magnitudes, or angle differences,  and/or (iv) reduce the approximation error based on iterative procedures \cite{Castillo2016}. More recently, these approximations can be applied to exploit the accuracy and inference speed of learning-based approaches, using the solution of historical (or synthetically generated) instances to predict base points \cite{Pineda2024}. 

Various approximations of the power flow equations rely on a truncated first-order Taylor series expansion of the nonconvex terms arising from the product of complex voltages. The difference between linear approximations using Taylor series lies in the selection of the independent variables, voltage magnitudes $v$ and phase angles $\theta$, to perform the truncated expansion. Zhang \ea \cite{Zhang2013} derives a linear power flow model using the $v$ and $\theta$ as independent variables. Yang \ea \cite{Yang2017b} refines such linear model by using $v^2$ and $\theta$, and consider the impact of $v$ on network losses in \cite{Yang2018}. Fatemi \ea \cite{Fatemi2015} use $v^2$ and $v^2\theta$ as independent variables and introduce an adjustment factor to minimize estimation errors.  Li \ea in \cite{Li2018} propose a logarithmic transform of voltage magnitudes, preserving the impact of voltage magnitudes and network losses on the active and reactive load flows. Yang \ea \cite{Yang2019} derive a general formulation, where the power injections are expressed as functions of the independent variables, for selecting independent variables that render a linear power flow with reduced errors. Fan \ea \cite{Fan2021} formulate a model to find a linear model with reduced linearization error by selecting the optimal independent variables for each bus. 
Long \ea in \cite{Long2024} propose a clustering-based approach to select the best linear power flow model depending on the similarity of the linearization error features of different operating conditions.

The authors in \cite{ONeill2012} and \cite{Castillo2016} propose the IV-Flow model, a linear approximation of the power flow equations by formulating the active and reactive power flow equations in rectangular coordinates in terms of voltages and currents and linearizing them using a first-order Taylor series expansion. Such approximation, along with the linearization of other constraints of the AC-OPF problem, is iteratively solved until a convergence criterion is satisfied. Mhanna and Mancarella \cite{Mhanna2022} propose a sequential linear programming approach based on the construction of supporting hyperplanes and halfspaces and the exclusive solution of linear programming problems. Coffrin and Van Hentenryck \cite{Coffrin2014} propose the LPAC model, which relies on a convex approximation of the cosine terms in the power flow equations and a Taylor approximation of the remaining nonlinear terms. Bolognani and Dörfler \cite{Bolognani2015} propose an implicit linearization of the power flow manifold with a bound on the worst-case error.

On the other hand, first-order Taylor series expansions can also be applied to derive convex approximations, where all the convexities are preserved in the model while the nonconvexities are linearized. Jabr in \cite{Jabr2006} proposes a second-order conic (SOC) approximation based on his SOC relaxation, where the quadratic equality constraints are relaxed into a convex inequality and the nonconvex (arctangent) constraints are linearized. Coffrin \ea in \cite{Coffrin2012} derive a quadratic programming approximation based on a second-order Taylor series expansion of the sine and cosine functions, assuming a near-nominal voltage magnitude approximation and a small angle difference, and on the relaxation of the quadratic equality constraints into inequalities. Baradar and Hesamzadeh in \cite{Baradar2015} propose a conic approximation of the OPF problem, assuming near nominal voltage magnitudes and small angle differences, expressed in terms of voltage magnitudes, power injections, and power losses. Šepetanc and Pandžić in \cite{Sepentanc2021} propose a second-order Taylor approximation of the AC power flows in polar coordinates by distributing power losses to both branch ends, and they assess its accuracy in a unit commitment problem.

Merkli \ea \cite{Merkli2021} formulate the AC power flow problem as a nonconvex quadratically constrained linear problem aiming at minimizing the deviation of power generation from given setpoints. The problem is reformulated as a difference-of-convex functions programming problem and solved using a difference-of-convex functions algorithm. The main limitations of such work are threefold: (i) the computational complexity of splitting the indefinite matrices in the nonconvex quadratic constraints, which relies on their full eigendecomposition (which takes $\mathcal{O}(n^3)$ time), (ii) the matrices to be split depend on the network topology, which makes such decomposition impractical in scenarios where the topology varies as it requires a new decomposition each time the topology is altered, and (iii) the convexified constraints do not preserve the sparsity of the original ones and are generally dense, which is computationally challenging for optimization solvers in terms of internal operations and memory usage.   

To bridge this gap, this paper proposes a quadratically-constrained convex approximation of the AC-OPF problem (QCAC) that does not rely on typical engineering assumptions such as high reactance-to-resistance ratio, near-nominal voltage magnitudes, or small angle differences, while preserving the structural sparsity of the original nonconvex constraints. The proposed QCAC approximation is derived by lifting the AC-OPF in rectangular coordinates, which allows reformulating the problem as a difference-of-convex functions program without requiring matrix-based decompositions. The reformulated problem is then convexified around a given set of nodal voltages by linearizing only the concave components of the voltage bilinearities.

In contrast to existing difference-of-convex approaches that rely on eigendecomposition of network-dependent matrices and iterative DC algorithms~\cite{Merkli2021}, the proposed method operates directly on individual bilinear voltage products, avoids dense matrix decompositions, and yields convex constraints that preserve the sparsity of the original formulation. As a result, the QCAC approximation scales favorably to large-scale networks and remains applicable under topology changes without requiring recomputation of network-dependent decompositions. The proposed approximation shows two important features: (i) it retains all the convex information of the nonconvex voltage relationships and solely linearizes their concave term, unlike other linear approximations, and (ii) it does not make any additional assumptions in terms of angle differences, reactance-to-resistance ratio, or voltage magnitudes. The accuracy and computational performance of the proposed approximation are numerically tested using realistic test systems up to 30,000 buses and two proof-of-concept case studies in transmission and distribution grids.

The remainder of this paper is organized as follows: Section~\ref{sec:proposed} presents the formulation of the AC-OPF problem in rectangular coordinates and derives the proposed QCAC approximation. Section~\ref{sec:experiments} illustrates the numerical performance of the proposed approximation with respect to the  SOC and SDP outer approximations and a first-order Taylor series linearization. Section~\ref{sec:case_study} demonstrates a case study of potential further applications of the proposed model. 
Section~\ref{sec:conclusion} concludes the paper.

\section{Proposed Convex Approximation of the AC-OPF Problem} \label{sec:proposed}
Next, we present the nomenclature and general formulation of the AC-OPF problem in lifted rectangular coordinates and derive the proposed approximation.

\subsection{Nomenclature}
\vspace{-1em} 

\nomenclature[A]{$\mathcal{G}$}{Generating units.}
\nomenclature[A]{$\mathcal{N}$}{Nodes in the electric network.}
\nomenclature[A]{$\mathcal{E}$}{Set of branches.}
\nomenclature[A]{$\mathcal{E}_i$}{Set of branches connected to node $i$.}

\nomenclature[B]{$\overline{B}_i/\underline{B}_i$}{Maximum/Minimum shunt susceptance at node $i$.}
\nomenclature[B]{$B_{ij}^{\rm sh}$}{Shunt susceptance of transmission line between nodes $i$ and $j$.}
\nomenclature[B]{$G_i^{\rm sh}/B_i^{\rm sh}$}{Shunt conductance/susceptance at node $i$.}
\nomenclature[B]{$G_{ij}/B_{ij}$}{Series conductance/susceptance of the branch between nodes $i$ and $j$.}
\nomenclature[B]{$\overline{N}_i^{\rm sh}$}{Maximum step of shunt compensation at node $i$.}
\nomenclature[B]{$P_i^{\rm d}/Q_i^{\rm d}$}{Active/Reactive power demand at node $i$.}
\nomenclature[B]{$\overline{P}_k/\underline{P}_k$}{Maximum/Minimum active power output of unit $k$.}
\nomenclature[B]{$\overline{Q}_k/\underline{Q}_k$}{Maximum/Minimum reactive power output of unit $k$.}
\nomenclature[B]{$\overline{S}_{ij}$}{Maximum apparent power flow through the branch between nodes $i$ and $j$.}
\nomenclature[B]{$T_{ij}$}{Tap ratio of the transformer between nodes $i$ and $j$.}
\nomenclature[B]{$\overline{V}_i/\underline{V}_i$}{Maximum/Minimum voltage magnitude at node $i$.}

\nomenclature[C]{$b_i^{\rm sh}$}{Shunt susceptance at node $i$.}
\nomenclature[C]{$c_{ii}$}{Squared voltage magnitude at node $i$.}
\nomenclature[C]{$c_{ij},s_{ij}$}{Lifted voltage variables.}
\nomenclature[C]{$n_i^{\rm sh}$}{Integer step of shunt compensation at node $i$.}
\nomenclature[C]{$p_k^{\rm g}/q_k^{\rm g}$}{Active/Reactive power produced by unit $k$.}
\nomenclature[C]{$p_{ij}/q_{ij}$}{Active/Reactive power flow from node $i$ to node $j$.}
\nomenclature[C]{$p_{ji}/q_{ji}$}{Active/Reactive power flow from node $j$ to node $i$.}
\nomenclature[C]{$q_i^{\rm sh}$}{Reactive power produced by shunt compensation at node $i$.}
\nomenclature[C]{$v_i^{\rm re}/v_i^{\rm im}$}{Real/Imaginary part of voltage at node $i$.}

\nomenclature[D]{$C_k(\cdot)$}{Cost function of unit $k$.}

\vspace{-1em}
{\printnomenclature}

\subsection{Problem Formulation}
The AC-OPF problem in lifted rectangular coordinates can be formulated as follows:
\begin{subequations} \label{eq:1}
\begin{align}
& \min_\Xi \quad \sum_{k\in \mathcal{G}} C_k (p_k^{\rm g}) \label{eq:1a} \\
& \subjto \quad  \nonumber \\
& \; \sum_{k\in\mathcal{G}_i}p_k^{\rm g} -P_i^{\rm d} - G_{i}^{\rm sh}c_{ii} = \smashoperator{\sum_{(i,j)\in\mathcal{E}_i}}p_{ij} + \smashoperator{\sum_{(j,i)\in\mathcal{E}_i}}p_{ji},\, \forall i \in \mathcal{N},  \label{eq:1b} \\
&\; \sum_{k\in\mathcal{G}_i}q_k^{\rm g} -Q_i^{\rm d} + B_{i}^{\rm sh}c_{ii} = \smashoperator{\sum_{(i,j)\in\mathcal{E}_i}}q_{ij} +\smashoperator{\sum_{(j,i)\in\mathcal{E}_i}}q_{ji},\,\forall i \in \mathcal{N}, \label{eq:1c} \\
&\; p_{ij} = \frac{G_{ij}}{T_{ij}^2}c_{ii} - \frac{G_{ij}}{T_{ij}}c_{ij} + \frac{B_{ij}}{T_{ij}}s_{ij}, \, \forall (i,j) \in \mathcal{E},  \label{eq:1d} \\
&\; p_{ji} = G_{ij}c_{jj} - \frac{G_{ij}}{T_{ij}}c_{ij} - \frac{B_{ij}}{T_{ij}}s_{ij},\, \forall (i,j) \in \mathcal{E}, \label{eq:1e} \\
&\; q_{ij} = -\frac{B_{ij} + B_{ij}^{\rm sh}}{T_{ij}^2}c_{ii} + \frac{G_{ij}}{T_{ij}}s_{ij} + \frac{B_{ij}}{T_{ij}}c_{ij},\, \forall (i,j) \in \mathcal{E},  \label{eq:1f} \\
&\; q_{ji} = -\left(B_{ij} + B_{ij}^{\rm sh}\right)c_{jj} - \frac{G_{ij}}{T_{ij}}s_{ij} + \frac{B_{ij}}{T_{ij}}c_{ij},\, \forall (i,j) \in \mathcal{E},  \label{eq:1g} \\
&\; \left(p_{ij}\right)^2 + \left(q_{ij}\right)^2 \le \overline{S}_{ij}^2,\, \forall (i,j),(j,i) \in \mathcal{E},  \label{eq:1h} \\
&\; \underline{V}_i^2 \le c_{ii} \le \overline{V}_i^2 ,\, \forall i \in \mathcal{N},     \label{eq:1i} \\
&\; \underline{P}_k \le p_k^{\rm g} \le \overline{P}_k,\, \forall k \in \mathcal{G}, \label{eq:1j} \\
&\; \underline{Q}_k \le q_k^{\rm g} \le \overline{Q}_k,\, \forall k \in \mathcal{G}, \label{eq:1k} \\
&\; c_{ii} = \big(v_i^{\rm re}\big)^2+\big(v_i^{\rm im}\big)^2,\, \forall i \in \mathcal{N}, \label{eq:1l} \\
&\; c_{ij} = v_i^{\rm re}v_j^{\rm re}+v_i^{\rm im}v_j^{\rm im},\, \forall (i,j) \in \mathcal{E}, \label{eq:1m} \\
&\; s_{ij} = v_i^{\rm re}v_j^{\rm im}-v_j^{\rm re}v_i^{\rm im},\, \forall (i,j) \in \mathcal{E}, \label{eq:1n}
\end{align}
\end{subequations}
where the optimization variables are elements of the set
\[
\Xi = \{p_k^{\rm g},q_k^{\rm g},p_{ij}^{\rm fr},q_{ij}^{\rm fr},p_{ij}^{\rm to},q_{ij}^{\rm to},c_{ii},c_{ij},s_{ij},v_i^{\rm re},v_i^{\rm im}\}.
\]

The objective function \eqref{eq:1a} aims at minimizing the total production cost. Constraints \eqref{eq:1b} and \eqref{eq:1c} enforce the active and reactive power balance, respectively. Constraints \eqref{eq:1d}-\eqref{eq:1g} define the active and reactive power flows throughout the power grid. Constraints \eqref{eq:1h} enforce the thermal capacity of the branches. Constraints \eqref{eq:1i} bound the squared voltage magnitudes. Constraints \eqref{eq:1j}-\eqref{eq:1k} bound the generating units' active and reactive power output, respectively. Constraints \eqref{eq:1l}-\eqref{eq:1n} are the nonconvex relationships between the voltage variables, $v_i^{\rm re}$ and $v_i^{\rm im}$, and the lifted voltage variables, $c_{ii}$, $c_{ij}$ and $s_{ij}$. Note that the nonconvexities of the AC-OPF are enclosed solely in the nonconvex quadratic constraints \eqref{eq:1l}-\eqref{eq:1n}.

\subsection{The QCAC Approximation of the AC Optimal Power Flow}

The formulation of the AC-OPF problem \eqref{eq:1} is amenable to be formulated as a difference-of-convex-functions programming problem, which is then convexified using a first-order Taylor series approximation of the concave terms. This subsection details the procedure to derive the proposed QCAC approximation.

The nonconvexities of Problem \eqref{eq:1} are in the form of bivariate products, $z = xy$ with continuous variables $x$ and $y$, which satisfy the following equality
\begin{equation}
    xy = \frac{1}{4}\left(x+y\right)^2 - \frac{1}{4}\left(x-y \right)^2, \label{eq:identity}
\end{equation}
which corresponds to the difference of two convex functions. This equality makes it suitable for the AC-OPF problem to be formulated as a difference-of-convex functions optimization problem. Thus, the nonconvex equality constraints \eqref{eq:1l}-\eqref{eq:1n} can be expressed as inequality constraints and we can reformulate each bilinear/quadratic term as a difference of convex functions using the equality \eqref{eq:identity}, which renders the following problem: 
\begin{subequations} \label{eq:3}
\begin{flalign}
& \min_\Xi \quad \eqref{eq:1a} \nonumber \\
& \subjto \quad \nonumber \\
&\quad \eqref{eq:1b} - \eqref{eq:1k},  \nonumber \\
&\quad c_{ii} \ge (v_i^{\rm re})^2 + (v_i^{\rm im})^2, \; \forall i \in \mathcal{N},  \label{eq:ccp_a} \\ 
&\quad c_{ii} \le \left(v_i^{\rm re}\right)^2+\left(v_i^{\rm im}\right)^2,\,i \in \mathcal{N}, \label{eq:ccp_b} \\
&\quad (v_i^{\rm re}+v_j^{\rm re})^2+(v_i^{\rm im}+v_j^{\rm im})^2 - 4c_{ij} \le \nonumber \\
&\quad \quad  (v_i^{\rm re}-v_j^{\rm re})^2+(v_i^{\rm im}-v_j^{\rm im})^2, \quad { \forall (i,j) \in \mathcal{E},  } \label{eq:ccp_c}\\
&\quad (v_i^{\rm re}-v_j^{\rm re})^2+(v_i^{\rm im}-v_j^{\rm im})^2 + 4c_{ij} \le \nonumber \\
&\quad \quad (v_i^{\rm re}+v_j^{\rm re})^2+(v_i^{\rm im}+v_j^{\rm im})^2, \quad  { \forall (i,j) \in \mathcal{E},  } \label{eq:ccp_d} \\
&\quad (v_i^{\rm re}-v_j^{\rm im})^2+(v_j^{\rm re}+v_i^{\rm im})^2 + 4s_{ij} \le \nonumber \\
&\quad \quad (v_i^{\rm re}+v_j^{\rm im})^2+(v_j^{\rm re}-v_i^{\rm im})^2, \quad  {\forall (i,j) \in \mathcal{E},  } \label{eq:ccp_e}\\
&\quad (v_i^{\rm re}+v_j^{\rm im})^2+(v_j^{\rm re}-v_i^{\rm im})^2 - 4s_{ij} \le \nonumber \\
&\quad \quad (v_i^{\rm re}-v_j^{\rm im})^2+(v_j^{\rm re}+v_i^{\rm im})^2, \quad  {\forall (i,j) \in \mathcal{E}}.& \label{eq:ccp_f} 
\end{flalign}
\end{subequations} 

Note that Problem~\eqref{eq:3} is an exact reformulation of the AC-OPF problem, which remains nonconvex due to the convex terms in the right-hand side of inequalities~\eqref{eq:ccp_b}-\eqref{eq:ccp_f}, which render these constraints nonconvex. Constraints~\eqref{eq:ccp_b}--\eqref{eq:ccp_f} are written as \( f^+(x) \le f^-(x) \), where both \( f^+ \) and \( f^- \) are convex functions. For inequalities of this form, convexity is preserved if \( f^+(x) \) is convex and \( f^-(x) \) is concave. Specifically, if \( f^-(x) \) is linear or concave, the inequality defines a convex constraint, since convex functions minus concave functions lead to convex inequalities. However, if \( f^-(x) \) is convex, the inequality will not define a convex constraint, and the problem will remain nonconvex. The proposed reformulation admits various (piecewise-)linearizations, as the nonconvexities in the right-hand side can be expressed as separable functions. A drawback of piecewise linearization is the need for SOS2 constraints, making the problem mixed-integer. 

To preserve the continuous nature of the AC-OPF problem, we propose linearizing these right-hand side convex terms using a first-order Taylor series approximation around given nodal voltages, $V_i^{\rm re}$ and $V_i^{\rm im}$, effectively convexifying constraints~\eqref{eq:ccp_b}-\eqref{eq:ccp_f}. It is important to note that the quality of the proposed approximation heavily depends on the choice of the linearization point, as it is the only source of information used to convexify the problem. Specifically, the quality of the QCAC approximation is highly sensitive to the choice of base voltages, \( V_i^{\rm re} \) and \( V_i^{\rm im} \), around which the first-order Taylor series linearization is performed. For standalone applications of the QCAC model in real-time operations, a reliable linearization point can typically be obtained from the most recent solution of the state estimator. Alternatively, high-quality convexification points can be predicted using learning-based algorithms, leveraging the solutions of historical instances to predict optimal convexification points. However, the further the operating point moves from the linearization point, the larger the error in the approximation, particularly in extreme system conditions where higher-order nonlinearities become significant. To mitigate this, the choice of base voltages should be made carefully to reflect the operating conditions, ideally based on recent measurements, or predicted operating states of the system.

The proposed QCAC approximation can be modeled as follows:

\begin{subequations} \label{eq:qcac}
\begin{align}
& \min_\Xi \quad \sum_{k\in \mathcal{G}} C_k (p_k^{\rm g}) + \rho\Bigg[\sum_{i \in \mathcal{N}} \xi_i^{\rm c} + \smashoperator{\sum_{(i,j) \in \mathcal{E}}} \left(\xi_{ij}^{\rm c} + \xi_{ij}^{\rm s} \right) \Bigg] \label{eq:qcac_a} \\
& \subjto \quad \eqref{eq:1b} - \eqref{eq:1k}, \nonumber \\
&\; c_{ii} \ge (v_i^{\rm re})^2 + (v_i^{\rm im})^2, \; \forall i \in \mathcal{N},  \label{eq:qcac_b} \\
&\; c_{ii} \le 2\left(V_i^{\rm re} v_i^{\rm re}+V_i^{\rm im} v_i^{\rm im} \right)-  \left(\Big(V_i^{\rm re}\Big)^2+\Big(V_i^{\rm im}\Big)^2\right) + \xi_i^{\rm c}, \; \forall i \in \mathcal{N},  \label{eq:qcac_c} \\
&\; (v_i^{\rm re}+v_j^{\rm re})^2+(v_i^{\rm im}+v_j^{\rm im})^2 - 4c_{ij} \le \xi_{ij}^{\rm c} + \nonumber \\
&\; \quad 2(v_i^{\rm re}-v_j^{\rm re})(V_i^{\rm re}-V_j^{\rm re})+  2(v_i^{\rm im}-v_j^{\rm im})(V_i^{\rm im}-V_j^{\rm im}) -\nonumber \\
&\; \quad \left[(V_i^{\rm re}-V_j^{\rm re})^2+(V_i^{\rm im}-V_j^{\rm im})^2\right], \; \forall (i,j) \in \mathcal{E},  \label{eq:qcac_d} \\
&\; (v_i^{\rm re}-v_j^{\rm re})^2+(v_i^{\rm im}-v_j^{\rm im})^2 + 4c_{ij} \le \xi_{ij}^{\rm c} + \nonumber \\
\displaybreak
&\; \quad 2(v_i^{\rm re}+v_j^{\rm re})(V_i^{\rm re}+V_j^{\rm re})+  2(v_i^{\rm im}+v_j^{\rm im})(V_i^{\rm im}+V_j^{\rm im})-\nonumber \\
&\; \quad \left[(V_i^{\rm re}+V_j^{\rm re})^2+(V_i^{\rm im}+V_j^{\rm im})^2\right] ,  \; \forall (i,j) \in \mathcal{E},  \label{eq:qcac_e} \\
&\; (v_i^{\rm re}-v_j^{\rm im})^2+(v_j^{\rm re}+v_i^{\rm im})^2 + 4s_{ij} \le \xi_{ij}^{\rm s} +\nonumber \\
&\; \quad 2(v_i^{\rm re}+v_j^{\rm im})(V_i^{\rm re}+V_j^{\rm im})+   2(v_j^{\rm re}-v_i^{\rm im})(V_j^{\rm re}-V_i^{\rm im}) -\nonumber \\
&\; \quad \left[(V_i^{\rm re}+V_j^{\rm im})^2+(V_j^{\rm re}-V_i^{\rm im})^2\right], \; \forall (i,j) \in \mathcal{E},  \label{eq:qcac_f} \\
&\; (v_i^{\rm re}+v_j^{\rm im})^2+(v_j^{\rm re}-v_i^{\rm im})^2 - 4s_{ij} \le \xi_{ij}^{\rm s} +\nonumber \\
&\; \quad 2(v_i^{\rm re}-v_j^{\rm im})(V_i^{\rm re}-V_j^{\rm im})+  2(v_j^{\rm re}+v_i^{\rm im})(V_j^{\rm re}+V_i^{\rm im}) - \nonumber \\
&\; \quad \left[(V_i^{\rm re}-V_j^{\rm im})^2+(V_j^{\rm re}+V_i^{\rm im})^2\right], \; \forall (i,j) \in \mathcal{E},  \label{eq:qcac_g} \\
&\; \xi_{i}^{\rm c} \ge 0,\; \forall i \in \mathcal{N},  \label{eq:qcac_h} \\
&\; \xi_{ij}^{\rm c} \ge 0,\; \forall (i,j) \in \mathcal{E},  \label{eq:qcac_i} \\
&\; \xi_{ij}^{\rm s} \ge 0,\; \forall (i,j) \in \mathcal{E},  \label{eq:qcac_j}
\end{align}
\end{subequations}
where $\xi_{i}^{\rm c}$, $\xi_{ij}^{\rm c}$, and $\xi_{ij}^{\rm s}$ are nonnegative slack variables and $\rho$ denotes a penalty coefficient.

If the slack variables are set to zero and the convexification points, $V_i^{\rm re}$ and $V_i^{\rm im}$, render a nonempty constraint set, Problem \eqref{eq:qcac} corresponds to an inner approximation of the AC-OPF problem. However, when the base points $V_i^{\rm re}$ and $V_i^{\rm im}$ render an empty constraint set, e.g., they are infeasible to the original constraint set, we can add the aforementioned nonnegative slack variables to relax such constraints. Notice that by adding the slack variables and penalizing them in the objective function, the problem is neither an inner nor an outer approximation of the original problem. Hence, it cannot provide valid upper or lower bounds on the optimal value of the objective function, respectively.

The proposed QCAC approximation can be also implemented using a sequential convex optimization procedure, summarized in Algorithm~\ref{alg:qcac}, which iteratively updates the convexification point and solves a sequence of convex QCQP subproblems.

\begin{algorithm}[h]
\caption{Sequential QCAC approximation for AC-OPF}
\label{alg:qcac}
\begin{algorithmic}[1]
\Require Base voltages $(V^{\rm re,0}, V^{\rm im,0})$, penalty coefficient $\rho^0>0$, maximum penalty coefficient $\rho^{\max}$, tolerance $\varepsilon>0$, and $\mu>1$
\Ensure Approximate AC-OPF solution $(v^{\rm re}, v^{\rm im}, p^{\rm g}, q^{\rm g})$

\State Set $t \gets 0$

\Repeat
    \State Form the QCAC subproblem \eqref{eq:qcac} by linearizing around $(V^{\rm re,t}, V^{\rm im,t})$ with penalty coefficient $\rho^{\rm t}$
    \State Solve the convex QCQP \eqref{eq:qcac} to obtain
    \[
        \Xi^{t+1} = \big(v^{\rm re,t+1}, v^{\rm im,t+1}, p^{\rm g,t+1}, q^{\rm g,t+1}, \xi^{\rm c,t+1}, \xi^{\rm c,t+1}_{ij}, \xi^{\rm s,t+1}_{ij}\big)
    \]
    \State Update base voltages:
    \[
        V^{\rm re,t+1} \gets v^{\rm re,t+1}, \qquad V^{\rm im,t+1} \gets v^{\rm im,t+1}
    \]
    \State Update penalty coefficient:
    $
    \rho^{\rm t+1} = \min \{\mu\rho^{\rm t}, \rho^{\rm max}\}
    $
    \State Check stopping criterion:
    $\sum_{i \in \mathcal{N}} \xi_i^{\rm c, t+1} + \smashoperator{\sum_{(i,j) \in \mathcal{E}}} \left(\xi_{ij}^{\rm c, t+1} + \xi_{ij}^{\rm s, t+1} \right) \le \varepsilon$
    \State $t \gets t+1$
\Until{stopping criterion is satisfied}

\State \Return $\Xi^{t} = \big(v^{\rm re,t}, v^{\rm im,t}, p^{\rm g,t}, q^{\rm g,t}, \xi^{\rm c,t}, \xi^{\rm c,t}_{ij}, \xi^{\rm s,t}_{ij}\big)$
\end{algorithmic}
\end{algorithm}

\subsection{Advantages of the QCAC approximation}
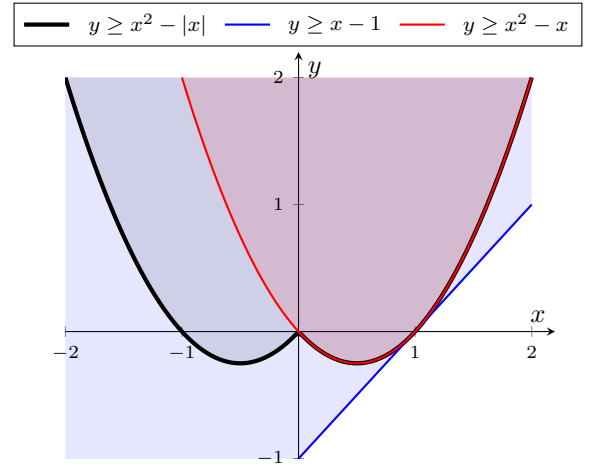
\begin{figure}[h]
\centering
\begin{tikzpicture}[scale=0.95]
    \begin{axis}[
        label style={font=\scriptsize},
        ticklabel style = {font=\scriptsize},
        axis lines = middle,
        xlabel = $x$,
        ylabel = $y$,
        ymin=-1, ymax=2.2,
        xmin=-2, xmax=2.2,
        samples=100,
        domain=-2:2,
        legend style={
            at={(0.5,1.1)},
            anchor=south,
            legend columns = 3,
            column sep=1ex,
            font = \scriptsize,
            align=left
        }
    ]
    
    \addplot [name path=A, domain=-2:2, black, ultra thick] {x^2 - abs(x)};
    \addlegendentry{$y \geq x^2 - |x|$}

    \addplot [name path=B, domain=-2:2, blue, thick] {(x-1)};
    \addlegendentry{$y \geq x-1$}
    
    \addplot [name path=C, domain=-1:2, red, thick] {x^2 - x};
    \addlegendentry{$y \geq x^2 - x$}

    \path[name path=top] (axis cs:-2,2) -- (axis cs:2,2);
    \path[name path=top2] (axis cs:-1,2) -- (axis cs:2,2);
    Shade regions
    \addplot [
        black,
        fill=black,
        fill opacity=0.1] 
    fill between[
        of=A and top,
        soft clip={domain=2:2}];

    \addplot [
        blue,
        fill=blue,
        fill opacity=0.1] 
        fill between[
            of=B and top,
            soft clip={domain=-2:2}];

    \addplot [
        red, 
        fill=red,
        fill opacity=0.1] 
        fill between[
            of=C and top2,
            soft clip={domain=-2:2}];
    \end{axis}
\end{tikzpicture}
\caption{Illustration of difference-of-convex function (nonconvex) constraint (black), and first-order Taylor series linearization (blue) and proposed convex approximation (red) around $x=1$. Note that the proposed convexification is tight for $x \ge 0$ whereas the first-order Taylor series linearization is tight only at the linearization point $x=1$.}
\label{fig:approximations}
\end{figure}

The main advantages of the proposed approximation over other linearization-based approaches, e.g., the sequential quadratic programming algorithm, and models based on linearization of the nonconvexities via first-order Taylor series expansion, are the following:
\begin{enumerate}
    \item The convexified constraints retain more information with respect to other linear approximations since they preserve all the information of the convex component of each nonconvex term and only linearize its concave part \cite{Lipp2016}. Conversely, in the linearization-based approaches, the second-order information is the objective constraints and/or constraints is lost. Figure~\ref{fig:approximations} illustrates in a simple difference-of-convex (nonconvex) constraint, $y \ge x^2 - \lvert x \rvert$, the approximate convex set obtained using the proposed convex approximation and the one obtained a standard first-order Taylor series linearization.
    \item The proposed approximation preserves the structure and sparsity of the original AC power flow equations. Each convexified constraint, \eqref{eq:qcac_b}-\eqref{eq:qcac_g}, depends on variables pertaining to the corresponding bus or branch. This fact enables the application of the QCAC approximation in settings where such feature is imperative, such as decentralized operation of power systems.
    \item The QCAC approximation can be used along learning-based techniques that can leverage historical (or synthetically generated) instances to learn the mapping between nodal demands and nodal voltages, which correspond to our convexification points \cite{Pineda2024}. Notice that since the QCAC approximation is convex, it can be embedded within the training pipeline of a learning model as a differentiable layer, enabling the use of packages like DiffOpt.jl \cite{Bensancon2024} and CVXPYLayers \cite{Agrawal2019} to automate the computation of the gradients.
\end{enumerate}

\section{Numerical experiments}\label{sec:experiments}

In this section, we provide numerical experiments to evaluate the performance of the proposed QCAC approximation in terms of feasibility, optimality gap, and runtime. The data of the studied instances can be retrieved from the \textsc{PGLib-OPF} library \cite{pglib}. The numerical experiments have been implemented on a Macbook Pro laptop with an Apple M2 Pro processor and 32 GB of RAM under JuMP 1.10 with Mosek 10.1.2, Ipopt, and Gurobi 11.0.2 as the optimization solvers for the convex, nonconvex, and mixed-integer nonconvex optimization problems, respectively.

\subsection{Experimental setup}

The results of the proposed QCAC approximation are compared with those obtained using a first-order Taylor series linearization (TS) of constraints \eqref{eq:1l}-\eqref{eq:1n}, as well as with a second-order conic (SOC) and a semidefinite (SDP) relaxation  \cite{Jabr2006,Bai2008}, for selected systems with up to 1354 buses. For the QCAC and TS models, the linearization/convexification points are set to the solution of the base case. For larger systems of 6,495, 13,659, and 30,000 buses, we study the performance of the QCAC model with respect to the solution of the original problem. The nonconvex, SOC, and SDP models are solved in PowerModels.jl in JuMP \cite{Coffrin2018}. For each system, we generate 100 demand samples by multiplying the base case's active and reactive power demand by a random vector drawn from a normal distribution with mean $1$ and standard deviation $0.1$. All numerical experiments in this paper are performed on balanced network models. The proposed QCAC approximation is formulated in rectangular voltage coordinates and is not intrinsically limited to single-phase systems; extension to three-phase unbalanced distribution networks would require introducing phase-specific voltage variables and applying the same difference-of-convex decomposition and linearization procedure to the resulting phase-coupled bilinear terms.

\subsection{Metrics}

\paragraph{Optimality gap} We assess the optimality gap of the cost of the projected feasible dispatch of the four models, i.e., QCAC, TS, SOC, and SDP, with respect to the cost of the AC-OPF problem. Such feasible dispatch, which is determined by projecting the approximate dispatches onto the constraint set of the AC-OPF problem, can be obtained by solving the following optimization problem:

 \begin{equation} \label{eq:proj}
\boldsymbol{\bar{p}}^{\rm CVX} \in \arg\min_{\Xi} \; \lVert \hat{\boldsymbol{p}}^{\rm CVX} - \boldsymbol{p}^g \rVert_2^2 \quad \text{s.t.} \quad \eqref{eq:1b} - \eqref{eq:1n}
\end{equation}
where the superscript $\rm CVX$ represents the corresponding linear/convex approximation/relaxations, i.e., $\rm CVX= \{QCAC,TS,SOC,$ $\rm SDP\}$, and $\boldsymbol{\bar{p}}^{\rm CVX}$ denotes the Euclidean projection of $\hat{\boldsymbol{p}}^{\rm CVX}$ onto the constraint set defined by \eqref{eq:1b}-\eqref{eq:1n}. 

The above optimization problem projects the active power generation dispatch to the feasible set of the AC-OPF problem. Note that Problem \eqref{eq:proj} is a nonconvex problem, and we solve it using Ipopt, which is a local solver that produces feasible solutions but not necessarily the closest feasible projection. The optimality gap is determined as follows:
\begin{equation*} \displaystyle
\text{Optimality gap}\,(\%) = 100 \frac{  \Big\lvert \sum\limits_{k \in \mathcal{G}}C_k\left(\overline{p}_k^{\rm CVX}\right) - \sum\limits_{k \in \mathcal{G}}C_k\left(p_k^{\ast}\right) \Big\rvert}{\sum\limits_{k \in \mathcal{G}}C_k\left(p_k^{\ast}\right)},
\end{equation*}
where $p_k^{\ast}$ denotes the solution of the AC-OPF problem \eqref{eq:1}.

\paragraph{Distance to feasibility} To measure how close the dispatch of the generators of the four models are with respect to the nonconvex constraint set of the AC-OPF problem, we define the distance to feasibility, which corresponds to the root mean square error between $\hat{\boldsymbol{p}}^{\rm g}$ and its orthogonal projection $\overline{\boldsymbol{p}}^{\rm g}$, with as follows:
\[
\text{Distance to feasibility}\,(\textrm{p.u.}) = \sqrt{\frac{1}{\lvert \mathcal{G} \rvert} \sum_{k\in \mathcal{G}}\left(\overline{p}_k^{\rm CVX} - \hat{p}_k^{\rm CVX} \right)^2}.
\]

\paragraph{Runtime} For the four models, we report the time to solve the relaxations and approximations in addition to the time to find the closest feasible dispatch, i.e., time to solve \eqref{eq:proj}.

\subsection{Results for small- and medium-scale instances} 

Table \ref{tab:results} presents the results for the small to medium-scale instances using the solution of the base case as a convexification point. The SDP relaxation is solved for instances up to 118 buses since the solution of larger instances violates the time limit of 600 seconds. The results show that the proposed approximation is competitive in terms of the optimality gap with respect to the other approximate models, particularly in the larger instances. In terms of distance to feasibility, the QCAC model notably outperforms the other three models in the majority of test cases, showing its ability to render a more secure generation dispatch. Notice that the quality of the solutions of the proposed approximation might (and often does) depend on how far the operating point is with respect to the convexification point.

\begin{table}[H]
\centering
\caption{Optimality Gaps, Distances to Feasibility, and Runtimes for Selected PGLIB Test Cases}
\renewcommand{\arraystretch}{1.05}
\label{tab:results}
\begin{adjustbox}{width=\textwidth}
\begin{tabular}{|l||c|c|c|c||c|c|c|c||c|c|c|c|}
\hline
                                      & \multicolumn{4}{c||}{Optimality Gap (\%)}                                          & \multicolumn{4}{c||}{Distance to feasibility  (p.u.)}                                & \multicolumn{4}{c|}{Runtime (s)}                                                                           \\
Test case                             & QCAC & TS & SOC & SDP & QCAC & TS & SOC & SDP & QCAC & TS & SOC & SDP \\
\hline \hline
pglib\_opf\_case30   & 0.06256            & 0.27451                    & 0.66478          & \textbf{0.00680}      
                     & 0.14951            & 0.39028                    & 0.29270          & \textbf{0.00004}      
                     & 0.03               & \textbf{0.03}              & 0.04             & 0.82   \\
pglib\_opf\_case39   & 0.27498            & 2.60725                    & 3.26863          & \textbf{0.21014}              
                     & \textbf{0.01517}   & 0.80445                    & 1.09835          & 0.07128               
                     & \textbf{0.03}      & \textbf{0.03}              & 0.04             & 2.82   \\
pglib\_opf\_case57   & 0.01706            & 0.05375                    & 0.00402          & \textbf{0.00092}      
                     & \textbf{0.01557}   & 0.61258                    & 0.01900          & 0.01852               
                     & 0.05               & \textbf{0.04}              & 0.04             & 12.39  \\
pglib\_opf\_case89   & 1.49676            &  7.79117                   & \textbf{1.14178} & 7.15318               
                     & \textbf{0.36341}   &  1.29853                   & 0.52210          & 1.18675          
                     & 0.37               & \textbf{0.35}              & 0.51             & 155.7  \\
pglib\_opf\_case118  & \textbf{0.37494}   & 10.01670                   & 2.89311          & 2.06640               
                     & \textbf{0.00300}   & 0.33364                    & 0.09849          & 0.06768               
                     & 0.18               & \textbf{0.17}              & 0.24             & 638.3  \\
pglib\_opf\_case179  & \textbf{1.21317}   & 0.04374                    & 2.46625          & - 
                     & \textbf{0.02196}   & 0.004736                   & 1.32447          & - 
                     & 0.58               & 0.36                       & \textbf{0.19}    & -     \\
pglib\_opf\_case300  & 1.53902            & 6.82574                    & \textbf{1.23162} & - 
                     & \textbf{0.18725}   & 0.39919                    & 0.48233          & - 
                     & \textbf{1.83}      & 2.21       & 2.16    & -    \\
pglib\_opf\_case500  & 0.31590            & 2.06582    & \textbf{0.13161} & - 
                     & 0.01643            & 0.14463    & \textbf{0.01076} & - 
                     & 1.05               & \textbf{1.00}       & 1.39             & -   \\
pglib\_opf\_case793  & \textbf{0.81594}   & 21.87617    & 3.85507          & - 
                     & \textbf{0.00300}   & 0.16165    & 0.09789          & - 
                     & 4.81      & \textbf{4.01}    & 4.32             & -      \\
pglib\_opf\_case1354 & \textbf{1.35456} &  57.87368    & 4.12804          & - 
                             & \textbf{0.01238}   &   1.61108     & 0.54144          & -
                             & \textbf{4.32}      &   4.32     & 5.09             & -     \\
\hline
\multicolumn{13}{l}{The highlighted values correspond to the best out of the QCAC, TS, SOC, and SDP models.} \\
\end{tabular}
\end{adjustbox}
\end{table}

\subsection{Results for large-scale instances}

We further assess the computational performance in terms of dispatch accuracy and solution speedup of the QCAC model with respect to its nonconvex counterpart using three large-scale test systems of up to 30,000 buses. 

\begin{figure}[ht]
    \centering
    \subfloat{%
       \includegraphics[width=0.475\columnwidth]{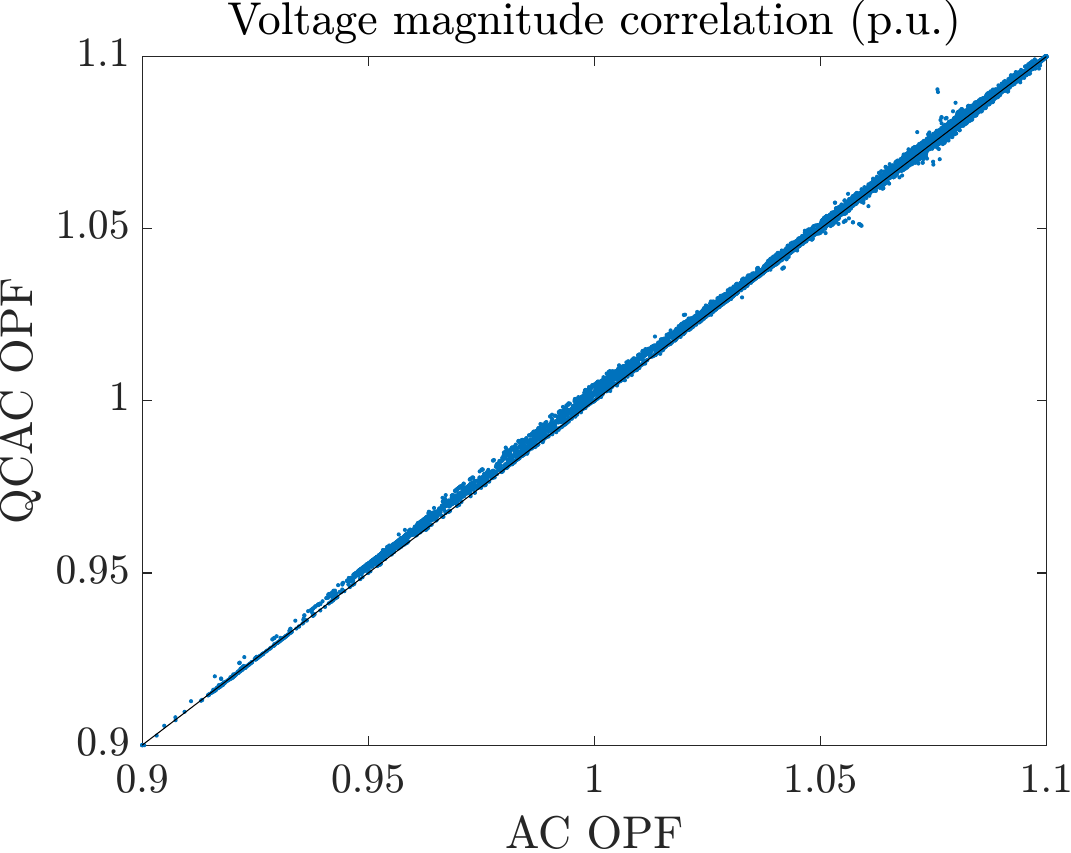}}
    \hfill
    \subfloat{%
       \includegraphics[width=0.47\columnwidth]{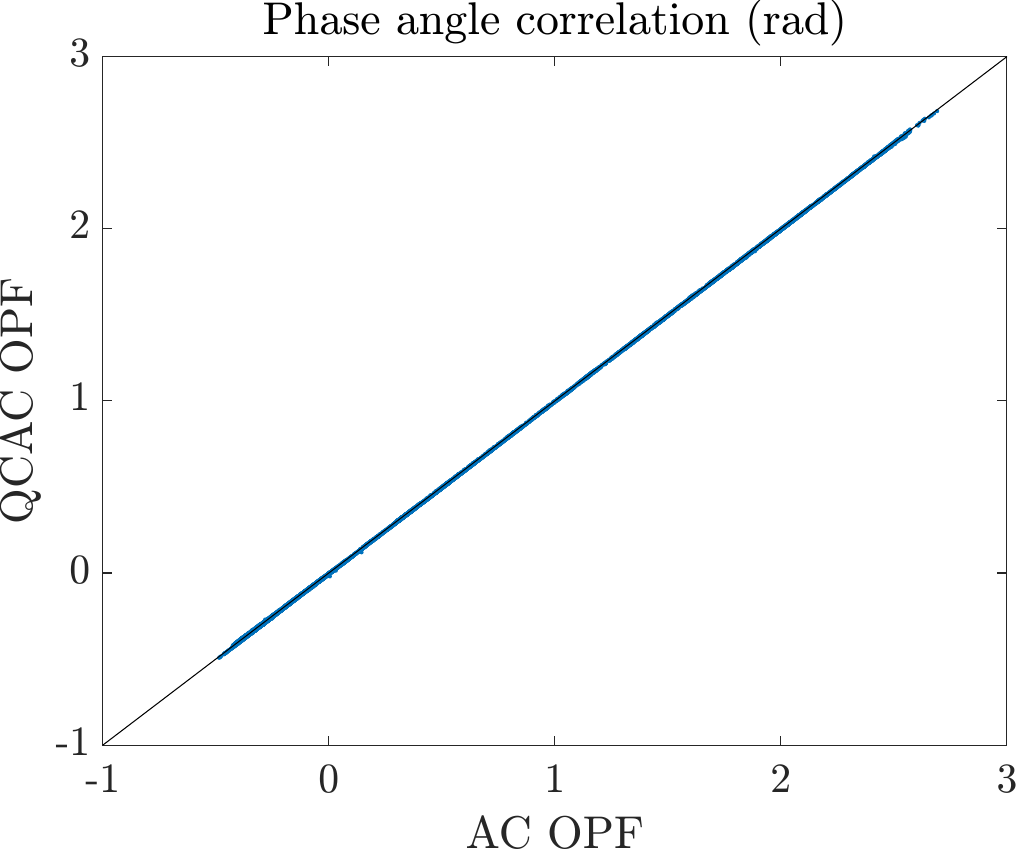}}
    \hfill
    \subfloat{%
       \includegraphics[width=0.475\columnwidth]{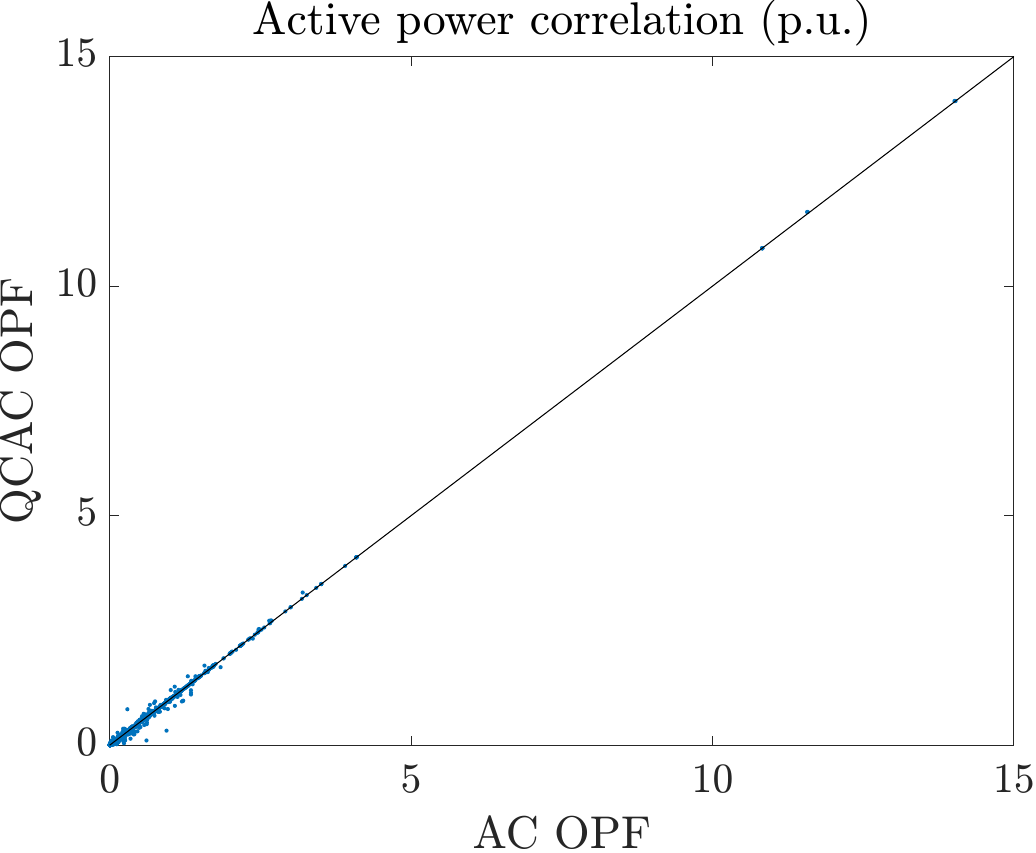}}
    \hfill
    \subfloat{%
       \includegraphics[width=0.47\columnwidth]{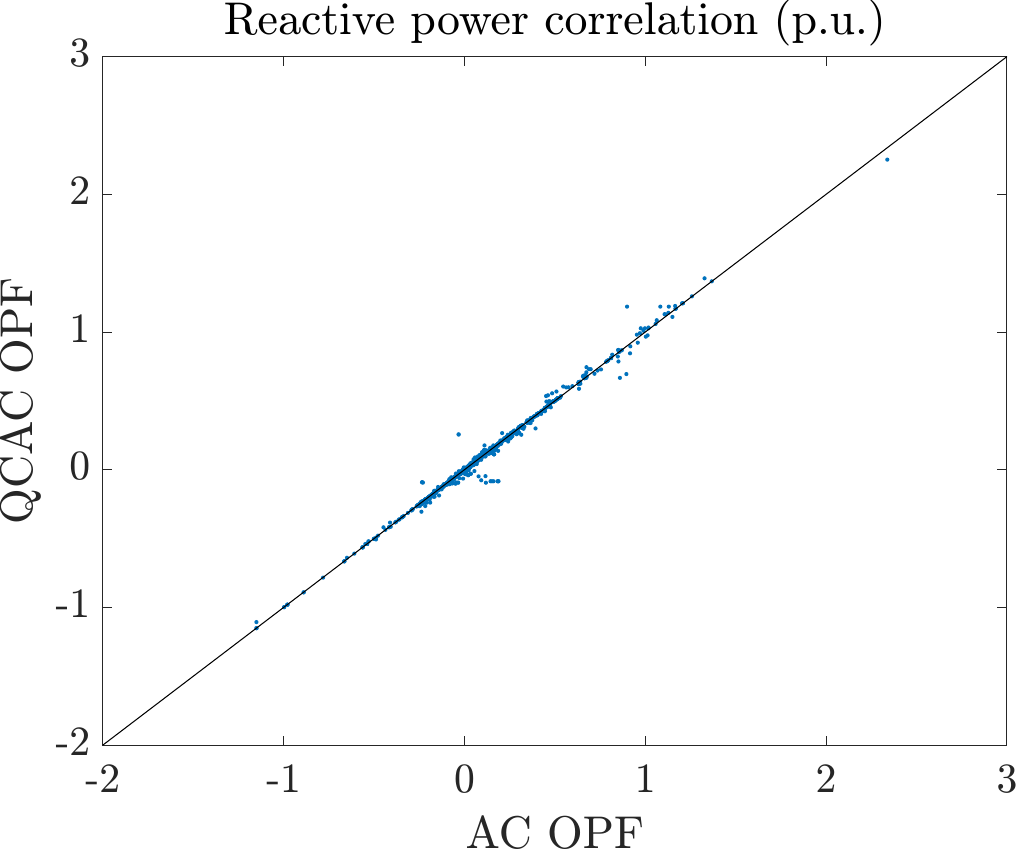}}
    \caption{Correlation of voltage magnitudes, phase angles, and active and reactive power generation for a randomly selected instance of the 30000-bus system}
    \label{fig:corr_large}
\end{figure}

Table~\ref{tab:error_stats_large} summarizes the correlation coefficients and mean/maximum errors of voltage magnitudes, phase angles, and active/reactive power generation. The QCAC model shows high accuracy across the three systems, showing its ability to properly convexify the nonlinearities of the AC-OPF problem. The reactive power is the most challenging variable in terms of both correlation coefficients and absolute errors. Figure~\ref{fig:corr_large} illustrates the high correlation of the four variables for an instance of the largest test system.

\begin{table}[H]
\centering
\caption{Correlation Coefficients, and Mean and Maximum Errors for Selected PGLIB Large-Scale Test Cases\label{tab:error_stats_large}}
\renewcommand{\arraystretch}{1.05}
\scriptsize
\begin{tabular}{|l||c|c|c|}
\hline
Test case & Corr. & Mean error & Max error \\
\hline \hline
\multicolumn{4}{|c|}{Voltage magnitude (p.u.)}                                                                                                                    \\
\hline
pglib\_opf\_case6495\_rte                  & 0.96352                                    & 0.00466                        & 0.22029                       \\
pglib\_opf\_case13659\_pegase              & 0.99348                                    & 0.00178                        & 0.15272                       \\
pglib\_opf\_case30000\_goc                 & 0.99915                                    & 0.00098                        & 0.03687                       \\
\hline 
\multicolumn{4}{|c|}{Phase angle (rad)}                                                                                                                          \\
\hline
pglib\_opf\_case6495\_rte                  & 0.98937                                    & 0.00466                        & 0.31219                       \\
pglib\_opf\_case13659\_pegase              & 0.99697                                    & 0.00178                        & 0.31959                       \\
pglib\_opf\_case30000\_goc                 & 0.99997                                    & 0.00098                        & 0.07397                       \\
\hline
\multicolumn{4}{|c|}{Active power (p.u.)}                                                                                                                         \\
\hline
pglib\_opf\_case6495\_rte                  & 0.99227                                    & 0.04500                        & 10.40000                      \\
pglib\_opf\_case13659\_pegase              & 0.99831                                    & 0.02053                        & 8.28580                       \\
pglib\_opf\_case30000\_goc                 & 0.99930                                    & 0.00665                        & 1.77600                       \\
\hline 
\multicolumn{4}{|c|}{Reactive power (p.u.)}                                                                                                                       \\
\hline
pglib\_opf\_case6495\_rte                  & 0.88023                                    & 0.06608                        & 7.99996                       \\
pglib\_opf\_case13659\_pegase              & 0.99338                                    & 0.01994                        & 3.19172                       \\
pglib\_opf\_case30000\_goc                 & 0.97869                                    & 0.00601                        & 2.92390                       \\
\hline
\end{tabular}
\end{table}

\begin{figure}[h]
    \centering
    \subfloat{%
       \includegraphics[width=0.7\columnwidth]{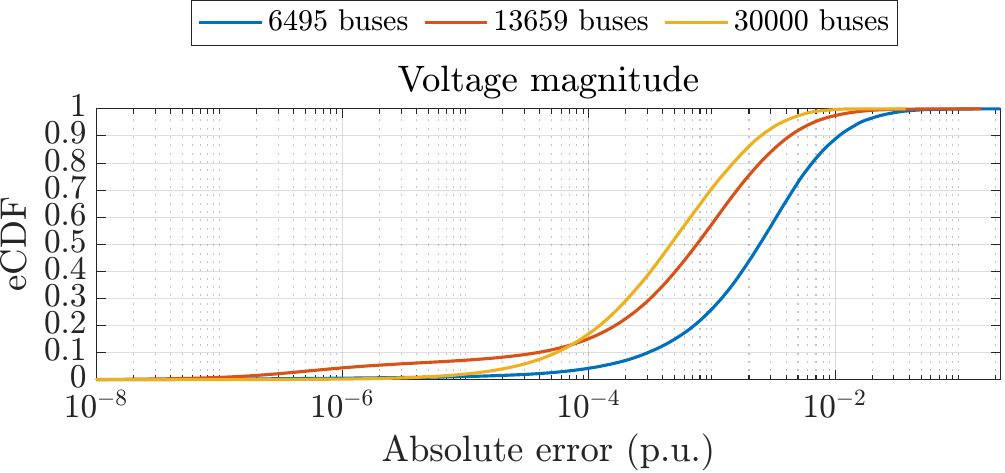}}
    \hfill
    \subfloat{%
       \includegraphics[width=0.7\columnwidth]{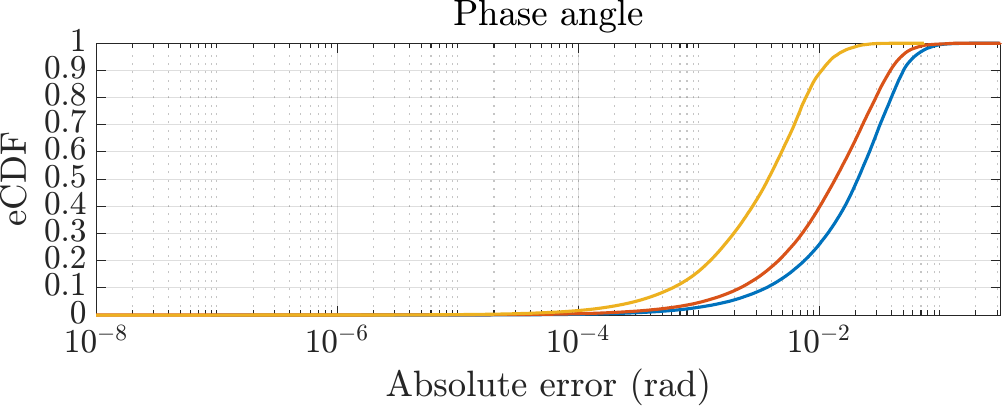}}
    \hfill
    \subfloat{%
       \includegraphics[width=0.7\columnwidth]{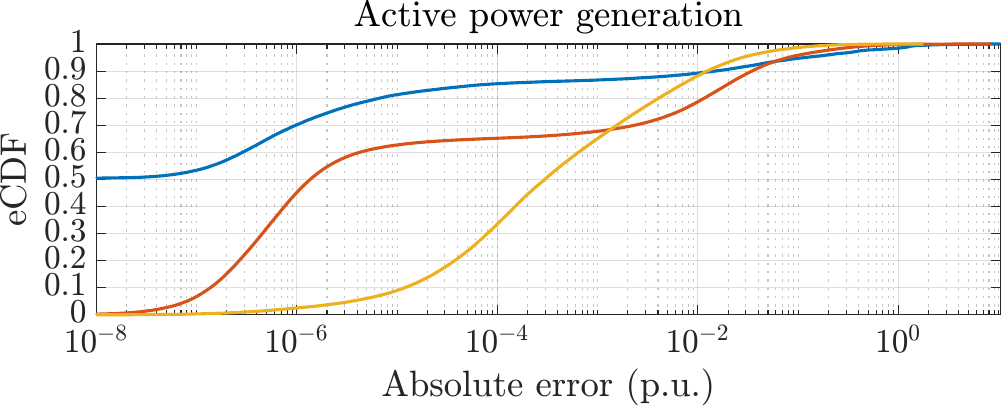}}
    \hfill
    \subfloat{%
       \includegraphics[width=0.7\columnwidth]{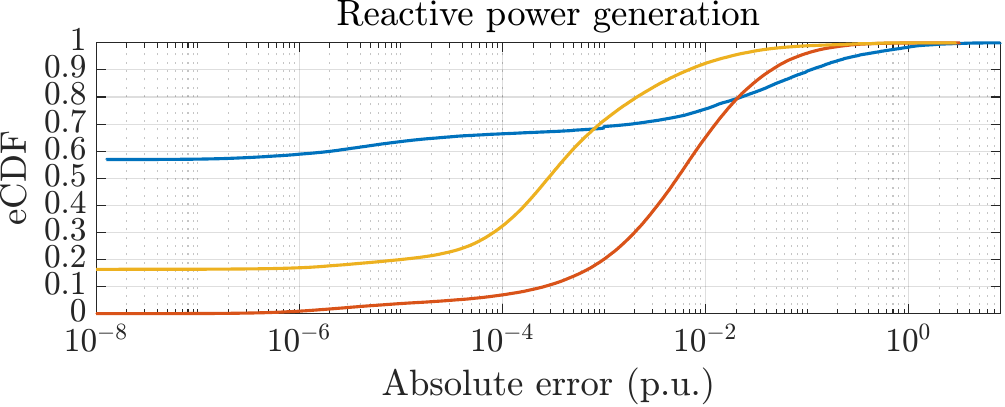}}
    \hfill
    \caption{Empirical cumulative distributions of absolute errors}
    \label{fig:ecdf_large}
\end{figure}

Figure~\ref{fig:ecdf_large} presents the empirical cumulative distribution functions for all the variables. For the active and reactive power generation of the largest system, more than $50\%$ of the instances show an absolute error of less than $1\times 10^{-8}$ p.u. The voltage magnitude errors for the three systems are below $0.011$ p.u. for 90\% of the instances.

\begin{figure}[h]
     \centering
     \includegraphics[width=0.6\columnwidth]{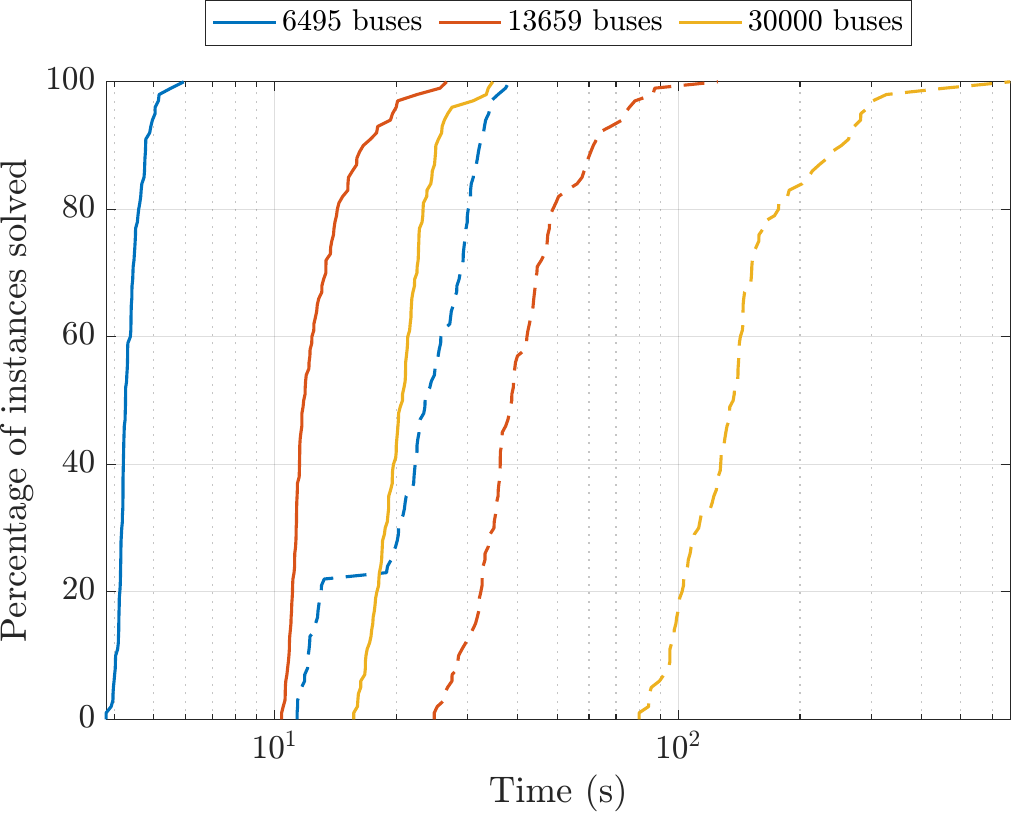}
     \caption{Time performance profile. QCAC OPF model (Solid). AC OPF model (Dashed).}
     \label{fig:time_large}
\end{figure}

Figure~\ref{fig:time_large} shows the time performance profile for the three systems for the QCAC model and the AC-OPF problem. The proposed model outperforms the nonconvex model in terms of solution time and reduced runtime variability. The largest system shows the highest runtime speedup, which empirically demonstrates the scalability of the proposed convexified model.

It is important to note that the quality of the solution of the proposed approximation relies heavily on the linearization point as it is the only source of information to convexify the problem. We note that for a standalone application of the QCAC model in real-time operations, a good linearization point can be easily obtained from the last solution of the state estimator. Another alternative to obtain high-quality convexification points is using end-to-end learning-based algorithms, which can predict such points by embedding the QCAC model into the learning pipeline and by leveraging the solution of historical instances.

\section{Case study} \label{sec:case_study}
In this section, we illustrate the application of the proposed approximation in the optimal reactive power flow with discrete controllers. The case study indicate that the proposed QCAC model provides an accurate convex approximation of the power flow equations, which is particularly useful when other approximations/relaxations are slow, inaccurate, or neglect voltage- and reactive-power-related variables and constraints. The numerical experiments were implemented considering the same setup as in Section~\ref{sec:experiments}.


The optimal reactive power flow (ORPF) problem, a mixed-integer nonconvex problem that derives from the OPF problem, aims at minimizing the generation costs by adjusting the active and reactive power generation setpoints as well as by operating discrete controllers, e.g., switchable shunt reactive power compensation \cite{Yang2017,Constante2021}. Operating such discrete controllers aims at acting on the nodal voltages by adjusting the reactive power injections to reduce the generating costs, requiring a detailed representation of both variables, i.e., voltage magnitudes and reactive power flows. Such requirement hinders the use of linear approximations of the power flow equations, such as the DC power flow, that neglect the aforementioned variables.  

A simplified formulation of the ORPF problem in lifted rectangular coordinates is presented below:
\begin{subequations} \label{eq:orpd}
\begin{align}
& \min_\Xi \quad \sum_{k\in \mathcal{G}} C_k (p_k^{\rm g}) \label{eq:orpd_a} \\
& \subjto \quad \nonumber \\
&\; \eqref{eq:1b}, \eqref{eq:1d} - \eqref{eq:1n}, \nonumber \\
&\; \sum_{k\in\mathcal{G}_i}q_k^{\rm g} -Q_i^{\rm d} + q_{i}^{\rm sh} = \smashoperator{\sum_{(i,j)\in\mathcal{E}_i}}q_{ij} +\smashoperator{\sum_{(j,i)\in\mathcal{E}_i}}q_{ji},\,\forall i \in \mathcal{N}, \label{eq:orpd_b} \\
&\; q_{i}^{\rm sh} = b_{i}^{\rm sh}c_{ii}, \, \forall i \in \mathcal{N}, \label{eq:orpd_c} \\
&\; b_{i}^{\rm sh} = \underline{b}_{i}^{\rm sh} + \left(\overline{b}_{i}^{\rm sh}-\underline{b}_{i}^{\rm sh} \right)\frac{n_i^{\rm sh}}{\overline{n}_i^{\rm sh}}, \, \forall i \in \mathcal{N}, \label{eq:orpd_d} \\
&\; n_i^{\rm sh} = \{0,1,\dots,\overline{n}_i^{\rm sh}\}, \, \forall i \in \mathcal{N}, \label{eq:orpd_e}
\end{align}
\end{subequations}
where the optimization variables are elements of the set:
\begin{align*}
\Xi^{\rm ORPF} = \{&p_k^{\rm g},q_k^{\rm g},p_{ij}^{\rm fr},q_{ij}^{\rm fr},p_{ij}^{\rm to},q_{ij}^{\rm to},q_i^{\rm sh},b_i^{\rm sh},n_i^{\rm sh}, c_{ii},c_{ij},s_{ij},v_i^{\rm re},v_i^{\rm im}\}.
\end{align*}

The objective function \eqref{eq:orpd_a} aims at minimizing the total production cost. Constraints \eqref{eq:orpd_b} enforce the reactive power balance constraint. Constraints \eqref{eq:orpd_c} define the reactive power injected by the shunt compensation per node whereas constraints \eqref{eq:orpd_d} define the shunt susceptance per node. Constraints \eqref{eq:orpd_e} model the discrete nature of the steps of shunt compensation per node. 

The nonconvexities of Problem \eqref{eq:orpd} arise from (i) the bilinear equalities, \eqref{eq:orpd_c}, (ii) the bilinear and quadratic equalities, \eqref{eq:1l}-\eqref{eq:1n}, and (iii) the discrete variables representing the steps of shunt compensation \eqref{eq:orpd_e}. The first sources of nonconvexities can be exactly linearized since $b_i^{\rm sh}$ are discrete, and both variables in the bilinear terms are bounded \cite[\S~9.2]{Williams2013}. The interaction of binary variables with the aforementioned nonconvexities makes impractical the solution of the original model for medium or large-scale systems.

To efficiently (approximately) solve Problem \eqref{eq:orpd}, we use the QCAC approximation, replacing constraints \eqref{eq:1l}-\eqref{eq:1n} by \eqref{eq:qcac_b}-\eqref{eq:qcac_j}, rendering in a mixed-integer convex optimization problem, which can be solved by state-of-the-art commercial solvers or well-documented decomposition techniques that exploit convexity, e.g., Benders decomposition. Note as well that the slack variables, $\xi_i^{\rm c},\xi_{ij}^{\rm c},\xi_{ij}^{\rm s}$, need to be penalized in the objective function.

To find the convexification point, $V_i^{\rm re}$ and $V_i^{\rm im}$, we solve the continuous relaxation of Problem \eqref{eq:orpd}. The experimental results study the benefits and accuracy of the convexified problem with respect to the original mixed-integer nonconvex formulation and the one using a SOC relaxation of the power flow equations. We use a 30-bus test system \cite{Alsac1974} under a normal operating condition (TYP) to test the performance of the QCAC approximation.  

The solution procedure for using the proposed convex approximation is as follows:
\begin{enumerate}
    \item Set $V_i^{\rm re}$ and $V_i^{\rm im}$ to the solution of the continuous relaxation of the ORPD problem.
    \item Solve the ORPD problem using the convexified constraints \eqref{eq:qcac_b}-\eqref{eq:qcac_j}.
    \item Fix $n_i^{\rm sh}$, $b_i^{\rm sh}$, and $q_i^{\rm sh}$ to the solution in the previous step, and solve the restricted Problem \eqref{eq:orpd}. 
\end{enumerate}

In the above procedure, we solve two continuous nonconvex optimization problems to find a convexification point and find a feasible solution for a fixed (near-)optimal value of the discrete controllers, respectively, and a mixed-integer convex problem to find such values for the discrete controllers.

Table \ref{tab:orpd_results} presents the numerical results for the ORPD problem. Gurobi was unable to find a feasible solution to the original MINLP model within 2 hours. Hence, we present the solution of its continuous relaxation, which provides a valid lower bound on the cost. To compute the cost and savings of the QCAC and SOC models, we fix the discrete variables obtained by the corresponding mixed-integer convex approximations and solve the restricted Problem \eqref{eq:orpd}, which provides a feasible solution. Note that the nonconvex problems are solved with Knitro, which cannot guarantee global optimality. 

The proposed approximation outperforms the one obtained by using the MISOCP relaxation in terms of cost savings, reducing the production cost by 19.43\% compared to the 14.49\% of the SOC relaxation. Note as well that techniques such as Benders decomposition, outer approximation, and an active set strategy can be used to reduce the computational burden of solving the problem and enable the scalability \cite{Constante2024,Jabr2024} of the SOC- and QCAC-based ORPD problems.

\begin{table}[H]
\centering
\caption{Results Summary: Proof-of-concept ORPD Problem\label{tab:orpd_results}}
\renewcommand{\arraystretch}{1.25}
\begin{tabular}{l|cc}
\hline \hline
\multirow{2}{*}{Model} & Cost  &  Savings   \\
                       & (\$)  &  (\%)      \\
\hline
AC-OPF                 & 4996.21 & -       \\     
$\text{ORPD}^{\dagger}$                   & 3974.46 & 20.45        \\
$\text{ORPD with SOC relaxation}^{\ddagger}$  & 4272.48 & 14.49         \\
$\text{ORPD with QCAC approximation}^{\ddagger}$  & \textbf{4025.62} & \textbf{19.43}    \\
\hline \hline 
\multicolumn{3}{l}{\scriptsize $\dagger$ The presented results correspond to the continuous relaxation.} \\
\multicolumn{3}{l}{\scriptsize $\ddagger$ The presented results are obtained following Procedure 1)-3).}
\end{tabular}
\end{table}

\subsection{Photovoltaic Hosting Capacity Problem}
We illustrate the application of the QCAC approximation to a 533-node radial distribution network whose data can be retrieved from Matpower \cite{Zimmerman2011}. Particularly, we focus on an application where the SOC relaxation of the power flow equations renders highly inaccurate solutions \cite{Low2014}. One such application is the photovoltaic (PV) hosting capacity problem, which aims to determine the maximum total PV generation that can be accommodated by a distribution feeder without violating its operation limits \cite{Jha2022}.  The PV hosting capacity problem is formulated as follows:
\begin{subequations} \label{eq:pvhc}
\begin{align}
& \min_{\Xi^{\rm PVHC}} \quad -\sum_{i\in \mathcal{N}} p_i^{\rm pv} \label{eq:pvhost_a} \\
& \subjto \quad  \nonumber \\
&\; \eqref{eq:1c}-\eqref{eq:1n}, \nonumber \\
& \; \sum_{k\in\mathcal{G}_i}p_k^{\rm g} + p_i^{\rm pv}-P_i^{\rm d} - G_{i}^{\rm sh}c_{ii} = \smashoperator{\sum_{(i,j)\in\mathcal{E}_i}}p_{ij} + \smashoperator{\sum_{(j,i)\in\mathcal{E}_i}}p_{ji}, \forall i \in \mathcal{N},  \label{eq:pvhost_b} \\
& \; \underline{P}^{\rm pv} \le p_i^{\rm pv} \le \overline{P}^{\rm pv}, \forall i \in \mathcal{N},  \label{eq:pvhost_c} \\
& \; p_{k}^{\rm g} \ge 0, \forall k \in \mathcal{G}_0,  \label{eq:pvhost_d} 
\end{align}
\end{subequations}
where the optimization variables are elements of the set:
\begin{align*}
\Xi^{\rm PVHC} = \{&p_k^{\rm g},q_k^{\rm g},p_i^{\rm pv},p_{ij}^{\rm fr},q_{ij}^{\rm fr},p_{ij}^{\rm to},q_{ij}^{\rm to}, c_{ii},c_{ij},s_{ij},v_i^{\rm re},v_i^{\rm im}\}.
\end{align*}

The objective function \eqref{eq:pvhost_a} is to maximize the sum of total power injected by the PV generators connected to a set of predetermined nodes. Constraints~\eqref{eq:pvhost_b} correspond to the nodal balance equations. Constraints~\eqref{eq:pvhost_c} bound the PV power injection per node, whereas constraint~\eqref{eq:pvhost_d} ensures that there is no reverse flow to the substation at the root node, $\mathcal{G}_0$,  of the feeder.

The approximated PV hosting problem using the proposed QCAC approximation, can be formulated as follows:
\begin{subequations} \label{eq:pvhc_qcac}
\begin{align}
& \min_{\Xi^{\rm PVHC}} \quad -\sum_{i\in \mathcal{N}} p_i^{\rm pv} + \rho\Bigg[\sum_{i \in \mathcal{N}} \xi_i^{\rm c} + \smashoperator{\sum_{(i,j) \in \mathcal{E}}} \left(\xi_{ij}^{\rm c} + \xi_{ij}^{\rm s} \right) \Bigg] \label{eq:pvhost_qcac_a} \\
& \subjto \quad  \nonumber \\
&\; \eqref{eq:1c}-\eqref{eq:1k}, \eqref{eq:qcac_b}-\eqref{eq:qcac_j}, \eqref{eq:pvhost_b}-\eqref{eq:pvhost_d},\nonumber 
\end{align}
\end{subequations}
where the optimization variables are elements of the set:
\begin{align*}
\Xi^{\rm PVHC} = \{&p_k^{\rm g},q_k^{\rm g},p_i^{\rm pv},p_{ij}^{\rm fr},q_{ij}^{\rm fr},p_{ij}^{\rm to},q_{ij}^{\rm to}, c_{ii}, c_{ij},s_{ij},v_i^{\rm re},v_i^{\rm im},\xi_i^{\rm c},\xi_{ij}^{\rm c},\xi_{ij}^{\rm s}\}.
\end{align*}

\begin{figure}[H]
     \centering
     \includegraphics[width=0.7\columnwidth]{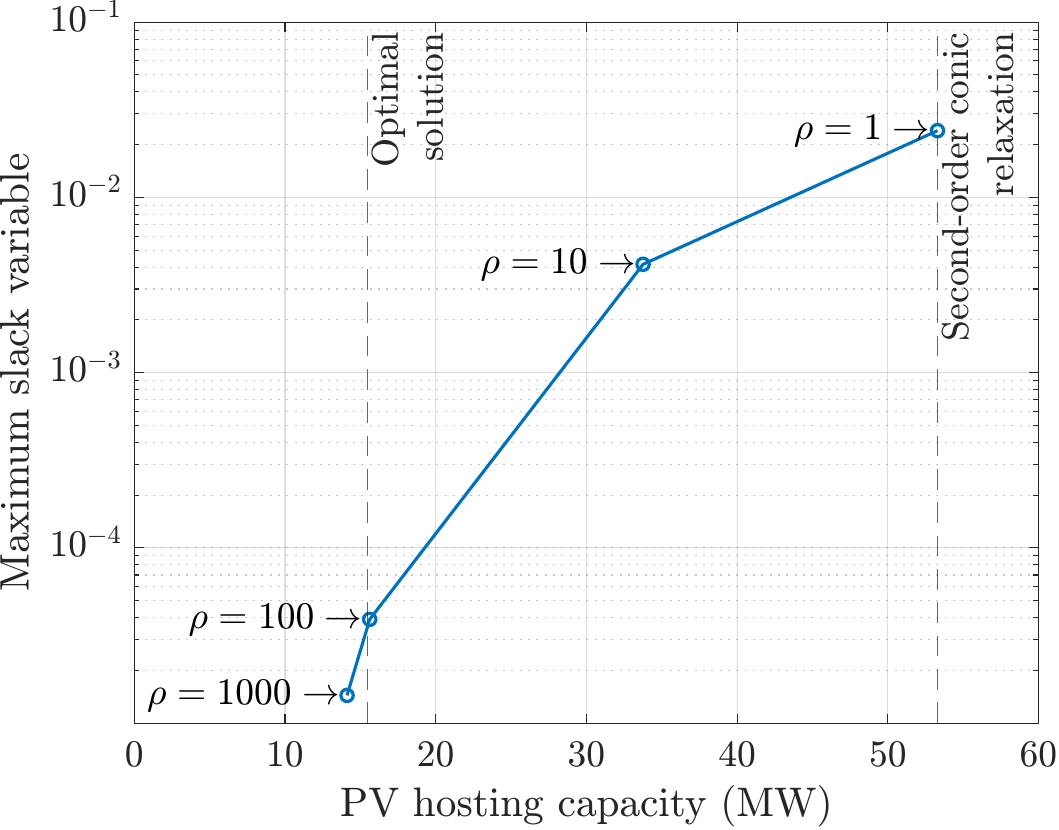}
     \caption{PV hosting capacity: Maximum slack-variable penalty trade-off curve.}
     \label{fig:pv_frontier}
\end{figure}

The convexification points of the power flow equations are set to a flat voltage profile. Figure \eqref{fig:pv_frontier} illustrates a frontier of the maximum slack variable and the PV hosting capacity for different values of the penalty coefficient $\rho$. Note that higher values of the penalty coefficient result in more conservative solutions that are closer to being feasible and to the convexification point. For $\rho=100$, the approximated problem results in a near-optimal and near-feasible solution with a maximum slack variable of $4\times10^{-5}$. For $\rho=1$, the approximated problem results in the same optimal objective function value as using the SOC relaxation with a maximum slack variable of $2.5\times10^{-2}$, which is highly infeasible. Note as well that according to the optimal value of the objective function, the SOC relaxation for this problem is a loose relaxation of the power flow equations.

Table~\ref{tab:pvhc_results} summarizes the results of the three models, i.e., nonconvex, SOC, and QCAC. For all the penalty coefficients, the proposed model offers a speedup of $\approx 33\times$ on average with the best cost error of $\approx 1\%$ whereas the SOC relaxation shows a higher speed up of $45\times$, but with a significantly higher cost error of $\approx 245\%$.

\begin{table}[H]
\centering
\caption{Result Summary: PV hosting capacity problem  \label{tab:pvhc_results}}
\label{tab:pvhc_summary}
\renewcommand{\arraystretch}{1.25}
\begin{tabular}{l|c|ccc}
\hline \hline
\multirow{2}{*}{Model}& \multirow{2}{*}{$\rho$} & Cost      & Cost error & Speedup \\
& & (MW) & (\%) & \\
\hline
Nonconvex             &        -            & 15.458          & -          & -       \\
\hline
SOC relaxation        &        -            & 53.300          & 244.80     & 44.8    \\
\hline
\multirow{4}{*}{QCAC approximation} & 1                   & 53.300          & 244.80     & 41.4    \\
                      & 10                  & 33.763          & 118.42     & 33.6    \\
                      & 100                 & \textbf{15.623} & \textbf{1.07}       & 29.5    \\
                      & 1000                & 14.118          & 8.67       & 27.5   \\
\hline \hline
\end{tabular}
\end{table}

\section{Conclusion}\label{sec:conclusion}

In this paper, we propose a quadratically-constrained convex approximation for the AC-OPF problem that relies solely on a convexification point, without the need for additional assumptions. The approach is based on a difference-of-convex reformulation of the power flow equations, followed by a convexification procedure using a first-order Taylor series expansion of the resulting concave terms. Our extensive numerical results demonstrate the advantages of this model, including a reduced optimality gap, near-feasible solutions, and smaller solution errors compared to the AC formulation. Furthermore, a case study illustrates the potential benefits of the proposed approximation over the second-order conic relaxation of the power flow equations for applications beyond the optimal power flow problem.

In our future work, we plan to make the proposed approximation adaptive in the sense that the convexification point and penalty coefficient are tuned to minimize the approximation error for a given operating range. We additionally plan to apply the proposed model to other applications such as an AC security-constrained unit commitment problem, transmission expansion planning problems, distribution OPF problem with distributed generation, and the robust OPF problem.

\section*{Declarations}

\textbf{Funding } 
This research was supported by the Davidson School of Chemical Engineering at Purdue University and the Office of Naval Research under grant no. N000142412641.

\noindent\textbf{Conflict of interest/Competing interests } 
The authors declare that they have no competing interests.

\noindent\textbf{Ethics approval and consent to participate } 
Not applicable.

\noindent\textbf{Consent for publication } 
Not applicable.

\noindent\textbf{Data availability } 
The data that support the findings of this study are available from the corresponding author upon reasonable request.

\noindent\textbf{Materials availability } 
Not applicable.

\noindent\textbf{Code availability } 
Not applicable.

\noindent\textbf{Author contributions } 
G.C.F. developed the proposed model, conducted the theoretical analysis, and implemented the numerical experiments. C.L. provided guidance on the methodology, supervised the research progress, and contributed to interpreting the results. G.C.F. and C.L. jointly wrote the manuscript and prepared all figures and tables. All authors reviewed and approved the final manuscript.

\bibliographystyle{spmpsci} 
{\small
\bibliography{bib/references.bib}}

\end{document}